\documentclass[dvips,psamsfonts,a4paper,11pt]{amsart}
\usepackage{amsfonts}%
\usepackage{amssymb}%
\usepackage{upref,xspace}
\usepackage[bookmarksopen, pdfstartview=FitH]{hyperref}
\theoremstyle{plain}
\newtheorem{thm}{Theorem}[section]
\newtheorem{cor}[thm]{Corollary}
\newtheorem{lem}[thm]{Lemma}
\newtheorem{prop}[thm]{Proposition}
\theoremstyle{definition}
\newtheorem{example}{Example}
\numberwithin{equation}{section}
\newcommand{\C}{\mathbb{C}}

\renewcommand{\P}{\mathbb{P}}
\renewcommand{\O}{\sheaf{O}}

\newcommand{\di}{\mathrm{d}}
\newcommand{\abs}[1]{\left|#1\right|}

\DeclareMathOperator{\diag}{diag}
\newcommand{\gr}[1]{\mathsf{#1}}
\newcommand{\alg}[1]{\mathfrak{#1}}
\newcommand{\term}[1]{\textit{#1}}
\newcommand{\sheaf}[1]{\mathcal{#1}}

\newcommand{\CR}{\ifmmode\mathrm{CR}\else\textrm{CR}\xspace\fi}
\makeatletter
\def\@strippedMR{}
\def\@scanforMR#1#2#3\endscan{%
  \ifx#1M\ifx#2R\def\@strippedMR{#3}%
  \else\def\@strippedMR{#1#2#3}%
  \fi\fi}
\renewcommand\MR[1]{\relax\ifhmode\unskip\spacefactor3000 
  \space\fi
  \@scanforMR#1\endscan
  \MRhref{\@strippedMR}{MR\,\@strippedMR}}
\makeatother
\renewcommand{\MRhref}[2]{%
  \href{http://www.ams.org/mathscinet-getitem?mr=#1}{#2}
}

\newcommand{\g}{\alg{g}}
\newcommand{\gz}{{\alg{g}_0}}
\newcommand{\q}{\alg{q}}
\newcommand{\gq}{(\gz,\q)}

\newcommand{\zo}{{(0,1)}}
\newcommand{\G}{\gr{G}}
\newcommand{\Q}{\gr{Q}}
\newcommand{\Gz}{\G_0}

\renewcommand{\o}{\mathbf{o}}
\begin{document}
%
\hypersetup{pdfauthor={Andrea Altomani}, 
pdftitle={Global CR functions on parabolic CR manifolds}}
\title[Global \CR functions on parabolic \CR manifolds]%
{Global \CR functions\\ on parabolic \CR manifolds}
\author{Andrea Altomani}
\address
{Dipartimento di Matematica\\
Universit\`a di Roma ``Tor Vergata''\\
Via della Ricerca Scientifica\\ 00133 Roma (Italy)}
\email{altomani@mat.uniroma2.it}%
\subjclass[2000]{Primary: 32V10; 
Secondary: 14M15, 32C09} 
\keywords{parabolic CR manifold, real group orbit, complex flag
  manifold, global CR functions}
\begin{abstract}Let $M$ be an orbit of a real semisimple Lie group $\Gz$
  acting on a complex a flag manifolds $\G/\Q$ of its complexification $\G$.
We study the space of global \CR functions on $M$ and characterize those
$M$ which are strictly locally \CR separable, i.e. those for which
global \CR functions induce local embeddings in $\C^n$.
\end{abstract}
\date{27 February 2007}
\maketitle
\section*{Introduction}\label{s:intro}%
Let $\G$ be a complex connected semisimple Lie group, $\Q$ a parabolic
subgroup and $\Gz$ a connected real form of $\G$.
The orbits of $\Gz$ in the flag manifold $Z=\G/\Q$ are finitely many,
and carry a natural \CR structure induced by their embedding in $Z$.
In \cite{AMN:pre} these $\Gz$-orbits are called \term{parabolic \CR 
manifolds}, and their \CR structure is invetigated by associating to them
their \CR algebras \cite{MN:2005}.

In this paper we describe the space of global \CR functions on parabolic
\CR manifolds.
This was done by Wolf \cite{Wolf:1969} in the case of open orbits; he
proved that an open orbit is holomorphically separable if and only if
it is a bounded symmetric domain.
For general orbits a similar statement is true, but there are some
important differences from the open case.

We need to introduce two notions of local \CR separability for a
\CR manifold $M$.
\term{Weak local \CR separability} means that global \CR functions locally
separate points, that is that there exist global \CR maps that give 
local \CR immersions into complex Euclidean spaces.
The stronger condition of  \term{strict local \CR separability} instead
requires the existence of global \CR maps into complex Euclidean spaces which
are locally  \emph{\CR embeddings}.
These two conditions are not equivalent, as Example~\ref{x:Athree} below
shows.

We prove that every parabolic \CR manifold $M$ admits a unique fibration 
$\rho\colon M\to M'$ over a parabolic \CR manifold $M'$ such that $M'$ is
strictly \CR separable and $\rho$ induces an isomorphism 
$\rho^*\colon\O(M')\to\O(M)$ of the spaces of
global \CR functions on $M$ and $M'$.

Then we characterize strictly locally \CR separable parabolic \CR manifolds
as those for which the fiber of the fundamental \CR
reduction is either a bounded symmetric domain or a smooth stratum in
the boundary of a bounded symmetric domain.

Finally we discuss some examples in detail.

I wish to thank prof. M. Nacinovich for several useful discussions
during the preparation of this paper

\section{Definitions and preliminary results}\label{s:CRalg}
Let $M=(M,HM,J)$ be a \CR manifold (see e.g. \cite{BER:1999} for the
relevant definitions).
We denote by $\O(M)$ the space of \emph{smooth \CR functions} on $M$.

We say that $M$ is (globally) \term{weakly \CR separable} if \CR
functions separate
points of $M$, that is if for every pair $x\neq y$ of points of $M$ 
there exists a \CR function $f\in\O(M)$ such that $f(x)\neq f(y)$.
We also say that $M$ is \term{weakly locally \CR separable} if every
point of $M$ has an open neighborhood $U$ such that global \CR functions on
$M$ separate points of $U$.

To introduce the notion of strict \CR separability of a \CR manifold
we need some more notation.
Let $\sheaf{S}=\sheaf{S}(M)\subset\Gamma(T^\C M)$ be the space of complex vector
fields $X$ on $M$ such that $X(f)=0$ for all $f\in\O(M)$, and let 
$S=S(M)\subset T^\C M$
be the vector distribution defined at $x\in M$ by:
\begin{equation}\label{e:defT}
S_x=\{X_x\mid X\in\sheaf{S}\}.
\end{equation}
By the definition of \CR functions, $T^\zo M\subset S$. 
We say that $M$ is \term{strictly locally \CR separable} at a point $x$
if $S_x=T_x^\zo M$, and that M is strictly locally \CR separable if it
is  strictly locally \CR separable at each point.
We have:

\begin{lem}\label{l:weakstrict}
  If $M$ is strictly locally \CR separable then $M$ is weakly locally \CR 
  separable.
\end{lem}
\begin{proof}
  Assume that there exist a point $p\in M$ and two sequences $x_n$ and
  $y_n$, with $x_n\neq y_n$ and converging to $p$, such that for every
  \CR function $f\in\O(M)$ we have $f(x_n)=f(y_n)$ for all $n$.
  Let $d$ be the distance function on $M$ defined by some Riemanniann
  metric $g$ on $M$.
  The functionals $\xi_n$, defined on a smooth function $f$ on $M$ by
  \begin{equation}
    \xi_n(f)=\left(f(x_n)-f(y_n)\right)/d(x_n,y_n),
  \end{equation} 
  converge,
  up to the choice of a subsequence, to a unit real tangent vector $X\in
  T_pM$.
  Clearly $X(f)=0$ for every $f\in\O(M)$, thus $M$ is
  not strictly locally \CR separable.
\end{proof}

Local strict \CR separability is an open condition, because $\dim S_x$
is upper semicontinuous with respect to $x$, and is actually equivalent
to the existence, for each point $x\in M$, of a
global \CR map of $M$ into a complex Euclidean space $\C^n$ that is a \CR
embedding in a neighborhood of $x$.

We turn now to the case of parabolic \CR manifolds.
We briefly recall the construction of a parabolic \CR manifold.
For a thorough study of these \CR manifolds we refer to \cite{AMN:pre}.

Let $\G$ be a connected complex semisimple Lie group, and $\Q$ a
parabolic subgroup of $\G$.
The homogeneous space $Z=\G/\Q$ is a compact simply connected smooth
complex projective variety, that is called a
\term{flag manifold}.
Fix a connected real form $\Gz$ of $\G$.
It acts on $Z$ by left multiplication, partitioning $Z$ into finitely
many $\Gz$-orbits.
A \term{parabolic \CR manifold} $M$ is one of these $\Gz$ orbits.

If $\g$, $\gz$ and $\q$ are the Lie algebras of $\G$, $\Gz$ and $\Q$
respectively, then the $\Gz$-orbit $\Gz\cdot\o$ through the point
$\o=e\Q$ of $Z$ is denoted by $M(\gz,\q)$ and the pair $\gq$ is the
\term{effective parabolic \CR algebra} associated to $M$ (see
\cite{MN:2005} for general definitions and properties of \CR algebras).
Since $\G$ acts transitively on $Z$ we can always, up to conjugation of
$\Q$ and $\Gz$ by an element of $\G$, assume that $\o\in M$. Hence every
$\Gz$-orbit in $Z$ is isomorphic to one of  the form $M=M(\gz,\q)$ for
some choice of $\gz$ and $\q$.

Since $\Gz$ acts on $Z$ by biholomorphisms, $M$ inherits from $Z$ a \CR
structure, and $\Gz$ acts on $M$ by real analytic \CR diffeomorphisms.
 
Let $\gq$ be an effective parabolic \CR algebra, and $M=M(\gz,\q)$ the 
associated parabolic \CR manifold.
Since $\sheaf{S}$ is invariant under \CR automorphisms of $M$, we have that
$S$ is a $\Gz$-homogeneous complex vector bundle on $M$.
Furthermore we have (see \cite[Lemma~14.5]{AMN:2006}):
\begin{lem}\label{l:subalgebra}
The vector subspace of $\g$:
\begin{equation}\label{e:defs}
	\alg{s}=(\di\pi^\C_e)^{-1}(S_{\o})
\end{equation}
 is a parabolic complex Lie subalgebra of
$\g$, containing $\q$.\qed
\end{lem}

Since $\q\subset\alg{s}$, we may consider the $\Gz$-equivariant fibration:
\begin{equation}\label{e:reduction}
	\rho\colon	M=M(\gz,\q)\to M'=M(\gz,\alg{s}).
\end{equation}
Every \CR function on $M'$ defines, via the pullback by $\rho$, a \CR function 
on $M$.
Indeed more is true:
\begin{thm}\label{t:fibration}
  Let $\gq$ be an effective parabolic \CR algebra and $M=M(\gz,\q)$ the 
  associated parabolic \CR manifold.
  Then there exists a unique $\Gz$-equivariant fibration
  $\rho\colon
  M\to M'$ onto a strictly locally \CR separable parabolic \CR
  manifold $M'$, such that $\rho$ induces an isomorphism
  on the space of \CR functions, that is:
  \begin{equation}
    \O(M)=\rho^{*}\O(M').
  \end{equation}

  The parabolic \CR manifold $M'$ is $\Gz$-equivariantly \CR isomorphic
  to $M(\gz,\alg{s})$, where $\alg{s}=(\di\pi^\C_e)^{-1}(S_{\o})$ is the
  Lie subalgebra of $\g$ defined in Lemma~\ref{l:subalgebra}.
\end{thm}
\begin{proof}
  From the definition of $\alg{s}$ it follows that 
  $M'=M(\gz,\alg{s})$ has all the
  required properties, so we only need to prove uniqueness.

  Let $\rho''\colon M\to M''=M(\gz,\q'')$ be a $\Gz$-equivariant
  fibration onto a strictly locally \CR separable parabolic \CR manifold
  $M''$ such that $\O(M)=\rho''^{*}\O(M'')$.
  Then $\q\subset\q''$ because $\rho''$ is a $\Gz$-equivariant
  fibration.
  On the other hand every \CR function on $M$ is the pullback of a \CR
  function on $M''$, which is  annihilated by $T^{0,1}M''$,
  hence $\q''\subset\alg{s}$.

  Finally $S(M)\subset\rho^*(S(M''))$ because $\rho$ is a \CR map and
  $S(M'')=T^{0,1}M''$ because $M''$ is strictly \CR separable.
  This yields $\q''=\alg{s}$, as stated.
\end{proof}
We refer to the $\Gz$-equivariant fibration \eqref{e:reduction}, or to
the corresponding $\gz$-equivariant fibration of \CR algebras, as the
\term{strictly \CR separable reduction} of $M(\gz,\q)$, or of $\gq$.

We can consider only simple parabolic \CR manifolds.
Indeed we have:
\begin{thm}\label{t:simple}
  Let $\gz=\bigoplus_i\gz_i$ 
  be the decomposition of the real semisimple Lie algebra $\gz$
  into the direct product of its simple ideals, and let $\g_i=\gz_i^{\C}$,
  $\q_i=\g_i\cap\q$.
  Then each $\q_i$ is parabolic in $\g_i$, $\q=\bigoplus_i\q_i$ and
  $M\gq$ is weakly (weakly locally, strictly locally) \CR separable
  if and only if all $M(\gz_i,\q_i)$'s are weakly (weakly locally, strictly
  locally) \CR separable.
\end{thm}
\begin{proof}
  The parabolic \CR manifold $M(\gz,\q)$ is isomorphic (see
  \cite[\S\,5.4]{Wolf:1969}) to the Cartesian product
  $\Pi_i M(\gz_i,\q_i)$.
\end{proof}

\section{Restriction to manifolds of finite type}
The following Theorem shows that we can restrict our consideration to
parabolic \CR manifolds of finite type.
\begin{thm}\label{t:fundamental}
  Let $M=M(\gz,\q)$ be a parabolic \CR manifold, 
  denote by $M'$ the fiber and by $M''=M(\gz,\q'')$ the
  base of a $\Gz$-equivariant \CR fibration $\phi\colon
  M\to M''$  onto a totally real 
  parabolic \CR manifold $M''$.
  Then $M$ is weakly (weakly locally, strictly locally) \CR separable if
  and only if $M'$ is weakly (weakly locally, strictly locally) \CR
  separable. 
\end{thm}
First we prove a lemma:
\begin{lem}\label{l:fundamental}
  Let $M=M(\gz,\q)$ be a parabolic \CR manifold, $M'$ the fiber and
  $M''=M(\gz,\q'')$ the base of a $\Gz$-equivariant \CR fibration 
  $\phi\colon
  M\to M''$  onto a totally real parabolic \CR manifold $M''$. 
  Then every $x\in M''$ has an open neighborhood $U\subset M''$ such
  that $\phi^{-1}(U)$ is \CR diffeomorphic by a real analytic map 
  to $U\times M'$.
\end{lem}
In particular Theorem~\ref{t:fundamental} and Lemma~\ref{l:fundamental}
apply when $\phi$ is the fundamental reduction \cite[\S\,5]{MN:2005} of $M$.
\begin{proof}[Proof of Lemma~\ref{l:fundamental}]
  Let $\alg{l}_0=\gz\cap\q$ be the isotropy subalgebra of $M$
  and $\alg{l}''_0=\gz\cap\q''$ that of $M''$.
  Since $M''$ is totally real, $\alg{l}''_0$ is a parabolic subalgebra of
  $\gz$, hence there exists a nilpotent subalgebra $\alg{n}_0$ complementary
  to $\alg{l}''_0$.
  Let $\gr{L}_0$, $\gr{L}''_0$, $\gr{N}_0$ be the analytic subgroups of
  $\gz$ with Lie algebras $\alg{l}_0$, $\alg{l}''_0$, $\alg{n}_0$,
  respectively and
  $\pi\colon\Gz\to M''=\Gz/\gr{L}_0''$ the projection onto the
  quotient. 
  
  The restriction of $\pi$ to $\gr{N}_0$ is a real analytic 
  local diffeomorphism.
  Choose an open neighborhood $W$ of the identity in $\gr{N}_0$ such that
  $\pi|_W$ is a diffeomorphism onto an open subset $\pi(W)=U$ of $M''$.
  Then the map:
  \begin{equation}
    \psi\colon U\times M'\ni (z,l\gr{L}_0)\mapsto 
    \left((\pi|_W)^{-1}(z)l\right)\gr{L}_0\in M
  \end{equation}
  is a real analytic \CR trivialization in a neighborhood of $e\gr{L}''$.
  The result follows because of the homogeneity of $M''$.
\end{proof}
\begin{proof}[Proof of Theorem~\ref{t:fundamental}]
  Let $x\neq y$ be two distinct points of $M$.
  If $\phi(x)\neq\phi(y)$ then we can choose any function $f$ on $M''$ 
  such that $f\circ\phi(x)\neq f\circ\phi(y)$, and $f\circ\phi$ is \CR,
  and separates $x$ and $y$.
  If $\phi(x)=\phi(y)$ then by Lemma~\ref{l:fundamental} we can find a
  \CR function $f$ on $M$ that separates $x$ and $y$ if and only if we
  can find such an $f$ on $\phi^{-1}(\phi(x))$.
  Thus $M$ is weakly (weakly locally) \CR separable if and only if 
  $M'$ is weakly (weakly locally) \CR separable.
  
  Fix a point $x\in M$, let $M'=\phi^{-1}(\phi(x))$  and denote by
  $\imath\colon M'\to M$ the inclusion map.
  Let $X\in T^{\C}_xM$ be a complex tangent vector at $x$ with
  $\di \phi^\C(X)\neq 0$ and $f$ a real analytic function on $M''$ such that
  $\di \phi^\C(X)(f)\neq 0$.
  Then $f\circ\phi$ is a \CR function on $M$ and $X(f\circ\phi)\neq 0$.
  This shows that $\sheaf{S}(M)=\imath_*(\sheaf{S}(M'))$, hence 
  $M$ is strictly locally \CR separable if and only if 
  $M'$ is strictly locally \CR separable.
\end{proof}

\section{Extension to Levi-flat orbits}
The case of totally complex parabolic \CR manifolds,
was discussed by Wolf in \cite{Wolf:1969}. 
There he proved (see \cite[Thm.~4.4.3]{FHW:2006}) the following:
\begin{prop}\label{p:open}
  Let $\gq$ be a simple totally complex parabolic effective \CR
  algebra and $M=M(\gz,\q)$ the corresponding totally complex
  parabolic \CR manifold.
  Then $M$ is weakly locally \CR separable if and only if $M$ is a
  bounded symmetric domain.
  In this case $M$ is also weakly \CR separable and strictly locally
  \CR separable.\qed
\end{prop}

We recall that a \CR manifold $M$ is \term{Levi-flat} if the analytic
tangent distribution $HM$ is integrable. For parabolic \CR manifolds this
is equivalent to the condition that the fibers of their fundamental
reduction are totally complex.
We have:
\begin{thm}\label{t:ext}
  Let $M=M(\gz,\q)$ be a $\Gz$-orbit in the complex flag
  manifold $Z=\G/\Q$.
  Then there exists a Levi-flat $\Gz$-orbit $N$ in $Z$, with
  $M\subset\overline{N}$, such that every \CR function $f$ on $M$
  continuously extends to a function $\tilde{f}$, continuous on $M\cup N$
  and \CR on $N$.

  If $M$ is of finite type, then $N$ is totally complex, hence open in $Z$.
\end{thm}
\begin{proof}
  If $M$ is Levi-flat we take $M'=M$.
  Otherwise, a theorem of Tumanov \cite{Tumanov:1990} asserts that there
  exists a complex wedge $W$, with edge contained in $M$, such that every \CR
  function $f$ on $M$ extends, continuously and uniquely, to a continuous 
  function $\hat{f}$ on $M\cup W$ that is holomorphic on $W$.
  Here by a complex wedge $W$ with edge in $M$ we mean a connected open
  subset $W$ of a complex submanifold $V$ of positive
  dimension of $Z$ such that $M\cap V$ is \CR generic in $V$ and 
  contained in the closure $\overline{W}$.

  Let $x\in W$ and define, for all $g\in\Gz$, a new $\tilde{f}$ by setting:
  \begin{equation}
    \tilde{f}(g\cdot x)=(\widehat{f\circ m_g})(x),
  \end{equation}
  where $m_g\colon M\to M$ denotes the action of $g$ on $M$.
  The function $\tilde{f}$ is well defined and \CR on the whole
  $\Gz$-orbit $M'=\Gz\cdot x$ through $x$. 
  By choosing $x$ close enough to $M$, we may arrange that
  $M\subset\overline{M'}$ and $\tilde{f}$ is continuous on $M\cup M'$. 

  By iterating this construction, we obtain a sequence of $\Gz$-orbits
  $M^{(i)}$ of nondecreasing dimension, each contained in the closure of
  the next.
  This sequence must necessary stabilize to a term $M^{(k)}=N$, that
  satisfies the first assertion of the theorem.

  If $M$ is of finite type, then also $N$ is of finite type and, being
  Levi flat and \CR generic, is open in $Z$.
\end{proof}

As a corollary, we obtain:
\begin{cor}\label{c:symmetric}
  If $M=M(\gz,\q)$ is a strictly locally \CR separable parabolic \CR
  manifold, embedded in the complex flag manifold $Z=\G/\Q$, then:
  \begin{enumerate}
  \item there exists a strictly locally \CR separable Levi-flat $\Gz$-orbit
    $N\subset Z$ with $M\subset\overline{N}$;
  \item $M$ is (globally) weakly \CR separable.
  \end{enumerate}
  If $M$ is of finite type then $N$ is a bounded symmetric domain.
\end{cor}
\begin{proof}
  We may assume that $M$ is the $\Gz$-orbit in $Z$
  through the point $\o=e\Q$.

  Let $N$ be the Levi-flat $\Gz$-orbit defined in Theorem~\ref{t:ext}.
  Let $\phi\colon Z\to
  Z'=\G/\Q'$ 
  be the $\G$-equivariant fibration of complex flag manifolds 
  that induces, by restriction to $N$, the strictly separable \CR
  reduction $\phi|_N\colon N\to N'\subset Z'$.
  Then every \CR function $f$ on $M$ extends continuously to a 
  function $\tilde{f}$, continuous on $M\cup N$ and \CR on $N$, constant
  along the fibers of $\phi$.
  By continuity also $f$ is constant along the fibers of $\phi$ and
  furthermore $f=\phi\circ f'$ for some \CR function $f'$ on $M'=\phi(M)$.
  This shows that:
  \begin{equation}
    S_{\o}(M)\supset (\di\phi^\C)^{-1}T^{(0,1)}_{\o}M'.
  \end{equation} 
  Since $M$ is strictly \CR separable, we obtain that:
  \begin{equation}
    T^{(0,1)}_{\o}M=(\di\phi^\C)^{-1}T^{(0,1)}_{\o}M',
  \end{equation} 
  which in turn implies that $\q=\q'$.
  Thus $N=N'$, that is $N$ is strictly locally \CR separable.

  If $M$ is of finite type, by Proposition~\ref{p:open}, $N$ is a
  bounded symmetric domain and $M\subset\overline{N}$.
  This fact also implies that $M$ is weakly \CR separable, thus the
  Theorem is proved if $M$ is of finite type.

  If $M$ is not of finite type, we apply two times
  Theorem~\ref{t:fundamental} to the fiber $M'$ of its fundamental
  reduction and obtain:
  \[\begin{split}
    \text{$M$ is strictly locally \CR separable} 
    \Longrightarrow\hspace{-5cm}&\\
    &\Longrightarrow\text{$M'$ is strictly locally \CR separable}
    \Longrightarrow\\
    &\Longrightarrow\text{$M'$ is weakly \CR separable}\Longrightarrow\\
    &\Longrightarrow\text{$M$ is weakly \CR separable,}
  \end{split}\]
  completing the proof.
\end{proof}

\section{Examples}
In this section we discuss some examples.
We will not need to utilize Tumanov's results, but more elementary
extension theorems will suffice.

In particular we recall the following statement.
Let $S^3=\{z\in\C^2\mid \abs{z}=1\}$ be the three-dimensional sphere,
endowed with the usual \CR structure, $B^2=\{z\in\C^2\mid\abs{z}<1\}$
the two-dimensional complex ball, 
$\Sigma$ a real two dimensional linear subspace of $\C^2$ (that may or
may not be a complex line) and set:
$\check{S}^3=S^3\setminus\Sigma$. 
Then every \CR function on $\check{S}^3$ extends continuously to a
function continuous on $\overline{B}^2\setminus\Sigma$ and holomorphic
on $\check{B}^2=B^2\setminus\Sigma$. 

For totally complex parabolic \CR manifolds and minimal parabolic \CR
manifolds the obstruction to \CR separability
is the esistence of embedded compact complex submanifold.
The general case is quite different.
Indeed we exhibit two examples of parabolic \CR manifold that are not
weakly locally \CR separable, but do not contain any compact complex
submanifold.

\begin{example}\label{x:Atwo} 
  Let $H(u,v)=u^*Av$ be the Hermitian form on $\C^3$ associated to the
  matrix $A=\diag(-1,1,1)$, and  $\G\simeq\gr{SL}(3,\C)$ the group of
  unimodular complex matrices. 
  The subgroup $\Gz$ of matrices in $\G$ that leave $H$ invariant is a
  real form of $\G$, isomorphic to $\gr{SU}(1,2)$.

  Let $Z$ be the complex flag manifold:
  \begin{equation}
    Z=\{\ell_1\subset\ell_2\subset\C^3\mid\dim\ell_i=i\}
  \end{equation}
  and $M$ the parabolic \CR manifold:
  \begin{equation}
    M=\{(\ell_1,\ell_2)\in Z\mid
    \text{$\ell_1$ is $H$-isotropic, $\ell_2$ is $H$-hyperbolic}\}.
  \end{equation} 
  Consider the $\G$-equivariant fibration 
  \begin{equation} 
    \phi\colon Z\ni (\ell_1,\ell_2)\mapsto \ell_1\in W
  \end{equation} 
  onto the complex flag manifold:
  \begin{equation}
    W=\{\ell_1\subset\C^3\mid\dim\ell_1=1\}.
  \end{equation}
  Then $\phi$ restricts to a $\Gz$-equivariant fibration $\phi\colon
  M\to N$ onto the 
  parabolic \CR manifold:
  \begin{equation}
    N=\{\ell_1\in W\mid\text{$\ell_1$ is $H$-isotropic}\}.
  \end{equation}
  The fiber of $\phi|_M$ over a point $\ell_1\in N$ is the set 
  $\{\ell_2\subset\C^3\mid
  \ell_1\subset\ell_2\not\subset(\ell_1)^\perp\}$, which is
  biholomorphic to $\C$.
  The \CR manifold $N$ is the boundary of the open domain:
  \begin{equation}
    D=\{\ell_1\in W\mid\text{$\ell_1$ is $H$-negative}\}.
  \end{equation}
  Fix an $H$-positive line $\hat{\ell}_1^+\subset\C^3$ and define:
  \begin{equation}
    M_{\hat{\ell}_1^+}=\{(\ell_1,\ell_2)\in
    M\mid\ell_2=\ell_1+\hat{\ell}_1^+\}
    \simeq\check{S}^3.
  \end{equation}
  Any \CR function $f$ on $M_{\hat{\ell}_1^+}$ extends continuously to a 
  function $\tilde{f}$, continuous on $M_{\hat{\ell}_1^+}\cup
  U_{\hat{\ell}_1^+}$ and holomorphic on 
  $U_{\hat{\ell}_1^+}$, where: 
  \begin{equation}
    U_{\hat{\ell}_1^+}=\{(\ell_1,\ell_2)\in Z\mid
    \text{$\ell_1$ is $H$-negative,
    $\ell_2=\ell_1+\hat{\ell}_1^+$}\}\simeq\check{B}^2.
  \end{equation}
  
  By letting $\hat{\ell}_1^+$ vary among all $H$-positive lines, we
  obtain that every \CR function $f$ on $M$ extends continuously to a 
  function $\tilde{f}$ continuous on $M\cup U$ and holomorphic on $U$, where:
  \begin{equation}
    U=\{(\ell_1,\ell_2)\in Z\mid
    \text{$\ell_1$ is $H$-negative and $\ell_2$ is $H$-hyperbolic}\}
  \end{equation}
  Let:
  \begin{equation}
  \begin{split}
    V&=\phi^{-1}(D)\setminus U\\&=\{(\ell_1,\ell_2)\in Z\mid
    \text{$\ell_1$ is $H$-negative, $H$ has signature $(0,{-})$ on $\ell_2$}\}.
  \end{split}
  \end{equation}
  Then $V$ has real codimension two in $\phi^{-1}(D)$ and is not complex
  analytic.
  Hence, by a theorem of Hartogs (\cite{Hartogs:1909}, see
  \cite[\S\,4, Thm.~3]{Narasimhan:1971}), there is a point $x\in V$ with
  the property that $\tilde{f}$ holomorphically extends to a full neighborhood
  $U_x$ of $x$ in $Z$.
  It follows that $\tilde{f}$ is constant on $\phi^{-1}\circ\phi(y)$ for
  all $y\in U_x$, and by  unique continuation $\tilde{f}$ is constant
  along the fibers of $\phi$, hence $M$ is not weakly locally \CR
  separable.
\end{example}
 
\begin{example}\label{x:Btwo} 
  Let $B(u,v)={}^tuAv$ and $H(u,v)=u^*Av$ be the bilinear and the Hermitian 
  form on $\C^5$, associated to the matrix $A=\diag(-1,-1,1,1,1)$, and 
  $\G\simeq\gr{SO}(5,\C)$ the group of unimodular complex matrices that
  preserve $B$.
  The connected component $\Gz$ of the identity in the subgroup of the real 
  matrices of $\G$ is isomorphic to $\gr{SO}^0(2,3)$ and the elements of $\Gz$
  also preserve the Hermitian form $H$.

  Let $Z$ be the complex flag manifold:
  \begin{equation}
    Z=\{\ell_1\subset\ell_2\subset\C^5\mid
    \text{$\dim\ell_i=i$, $\ell_i$ is $B$-isotropic}\}
  \end{equation}
  and $M$ the parabolic \CR manifold:
  \begin{equation}
    M=\{(\ell_1,\ell_2)\in Z\mid
    \text{$\ell_1$ is $H$-isotropic, $\ell_1\neq\bar\ell_1$, $\ell_2$ is 
    $H$-hyperbolic}\}.
  \end{equation} 
  Consider the $\G$-equivariant fibration 
  \begin{equation} 
    \phi\colon Z\ni (\ell_1,\ell_2)\mapsto \ell_1\in W
  \end{equation} 
  onto the complex flag manifold:
  \begin{equation}
    W=\{\ell_1\subset\C^5\mid\text{$\dim\ell_1=1$, $\ell_1$ is $B$-isotropic}\}.
  \end{equation}
  Then $\phi$ restricts to a $\Gz$-equivariant fibration $\phi\colon
  M\to N$ onto the 
  parabolic \CR manifold: 
  \begin{equation}
    N=\{\ell_1\in W\mid
    \text{$\ell_1$ is $H$-isotropic, $\ell_1\neq\bar\ell_1$}\}.
  \end{equation}
  The \CR manifold $N$ is an open stratum in the boundary of the open domain:
  \begin{equation}
    D=\{\ell_1\in W\mid\text{$\ell_1$ is $H$-negative}\}.
  \end{equation}
  Fix a $B$-isotropic and $H$-positive line $\hat{\ell}_1^+\subset\C^5$
  and define: 
  \begin{equation}
    M_{\hat{\ell}_1^+}=\{(\ell_1,\ell_2)\in
    M\mid\ell_2=\ell_1+\hat{\ell}_1^+\}
    \simeq\check{S}^3.
  \end{equation}
  Any \CR function on $M_{{\ell}_1^+}$ continuously extends to a
  function $\tilde{f}$, continuous on
  $M_{\hat{\ell}_1^+}\cup U_{\hat{\ell}_1^+}$ and hoomorphic on
  $U_{\hat{\ell}_1^+}$, where: 
  \begin{equation}
    U_{\hat{\ell}_1^+}=\{(\ell_1,\ell_2)\in Z\mid
    \text{$\ell_1$ is $H$-negative, $\ell_2=\ell_1+\hat{\ell}_1^+$}\}
    \simeq\check{B}^2.
  \end{equation}
  
  By letting $\hat{\ell}_1^+$ vary among all $B$-isotropic and
  $H$-positive lines, we obtain that every \CR function $f$ on $M$
  continuously extends to a function $\tilde{f}$, continuous on
  $M\cup
  U$ and holomorphic on $U$, where: 
  \begin{equation}
    U=\{(\ell_1,\ell_2)\in Z\mid
    \text{$\ell_1$ is $H$-negative and $\ell_2$ is $H$-hyperbolic}\}.
  \end{equation}

  Let:
  \begin{equation}
  \begin{split}
	    V&=\phi^{-1}(D)\setminus U\\&=\{(\ell_1,\ell_2)\in Z\mid
	    \text{$\ell_1$ is $H$-negative and $\ell_2$ is $H$-degenerate}\}.
  \end{split}
  \end{equation}
  Then $V$ has real codimension two in $\phi^{-1}(D)$ and is not complex
  analytic,
  hence by using the same argument of 
  Example~\ref{x:Atwo}, it follows that any \CR function on $M$ is
  constant along the fibers of $\phi$. Thus $M$ is not weakly locally \CR
  separable.
\end{example}
 
The next example consists of a parabolic \CR manifold that is weakly,
but not strictly, locally \CR separable.

\begin{example}\label{x:Athree}
  Let $\G\simeq\gr{SL}(4,\C)$ be the group of unimodular $4\times 4$ complex 
  matrices and $\Gz\simeq\gr{SU}(2,2)$ the subgroup of matrices 
  preserving the Hermitian form $H$ associated to the matrix 
  $\diag(-1,-1,1,1)$.

  Let $Z$ be the complex flag manifold:
  \begin{equation}
    Z=\{\ell_1\subset\ell_2\subset\C^5\mid\dim\ell_i=i\},
  \end{equation}
  and $M$ the parabolic \CR manifold:
  \begin{equation}
    M=\{(\ell_1,\ell_2)\in Z\mid\text{$\ell_1$ is $H$-isotropic, 
    $H|_{\ell_2}$ has signature $(0,{-})$}\}.
  \end{equation}
  Consider the $\G$-equivariant fibration 
  \begin{equation}
    \phi\colon Z\ni (\ell_1,\ell_2)\mapsto \ell_2\in W
  \end{equation} 
  onto the complex flag manifold:
  \begin{equation}
    W=\{\ell_2\subset\C^5\mid\dim\ell_2=2\}.
  \end{equation}
  Then $\phi$ restricts to a $\Gz$-equivariant fibration $\phi\colon M\to N$,
  which is a \CR map  and a smooth diffeomorphism, but not a \CR fibration, 
  onto the parabolic \CR manifold:
  \begin{equation}
    N=\{\ell_2\in W\mid\text{$H|_{\ell_2}$ has signature $(0,{-})$}\}.
  \end{equation}
  The \CR manifold $N$ is an open stratum in the boundary of the open domain:
  \begin{equation}
    D=\{\ell_2\in W\mid\text{$\ell_2$ is $H$-negative}\}.
  \end{equation}
  Fix an $H$-negative line $\hat{\ell}_1^-\subset\C^5$ and define:
  \begin{equation}
    M_{\hat{\ell}_1^-}=\{(\ell_1,\ell_2)\in M\mid 
    \text{$\ell_1$ is $H$-isotropic, $\ell_1\perp_H\hat\ell_1^-$,
    $\ell_2=\ell_1+\hat\ell_1^-$}\}\simeq{S}^3.
  \end{equation}
  Any \CR function on $M_{{\ell}_1^-}$ continuously extends to a
  function continuous on $M_{\hat{\ell}_1^-}\cup U_{\hat{\ell}_1^-}$,
  and holomorphic on $U_{\hat{\ell}_1^-}$, where:
  \begin{equation}
    U_{\hat{\ell}_1^-}=\{(\ell_1,\ell_2)\in Z\mid
    \text{$\ell_1$ is $H$-negative, $\ell_1\perp_H\hat\ell_1^-$,
    $\ell_2=\ell_1+\hat{\ell}_1^-$}\}\simeq{B}^2.
  \end{equation}

  By letting $\hat{\ell}_1^-$ vary among all $B$-isotropic and
  $H$-positive lines, we obtain that every \CR function on $M$
  continuously extends to a function continuous on $M\cup U$ and
  holomorphic on $U$, where: 
  \begin{equation}
    U=\{(\ell_1,\ell_2)\in Z\mid
    \text{$\ell_2$ is $H$-negative}\}.
  \end{equation}
  
  Since each fiber of the restriction of $\phi$ to $U$ is biholomorphic
  to $\C\P^1$, every \CR function $f$ on $M$ can be extended to a
  \CR function $\tilde{f}$ on $\phi^{-1}(N)$, which is constant along the
  fibers of $\phi$.
  This shows that $f$ is also \CR on $N$, hence $M$ is not strictly
  locally \CR separable.
  
  On the other hand $N$ is strictly locally \CR separable, hence
  by Lemma~\ref{l:weakstrict} $M$ is weakly locally \CR separable.
\end{example}

\bibliographystyle{amsplain}
\bibliography{crfunc}
\end{document}